\documentclass[11pt]{amsart}

\usepackage{fullpage}

\usepackage{amssymb,floatflt,amscd}

\newtheorem{theorem}[subsection]{Theorem}
\newtheorem{proposition}[subsection]{Proposition}
\newtheorem{lemma}[subsection]{Lemma}
\newtheorem{corollary}[subsection]{Corollary}
\newtheorem{definition}[subsection]{Definition}

\newcommand{\1}{{1\!\!1}}

\newcommand{\Z}{\mathbb{Z}}
\newcommand{\C}{\mathbb{C}}
\newcommand{\N}{\mathbb{N}}
\newcommand{\Q}{\mathbb{Q}}

\renewcommand{\o}{\otimes}

\DeclareMathOperator{\Aut}{Aut}

\DeclareMathOperator{\Map}{Map}

\DeclareMathOperator{\Ind}{Ind}
\DeclareMathOperator{\Res}{Res}
\DeclareMathOperator{\Spec}{Spec}

\newcommand{\CM}{\mathcal{M}}

\newcommand{\Mbar}{\overline{\mathcal{M}}}
\renewcommand{\]}{{]\!]}}

\newcommand{\bull}{\bullet}
\renewcommand{\SS}{\mathbb{S}}
\renewcommand{\P}{\mathbb{P}}

\newcommand{\A}{\mathbb{A}}
\DeclareMathOperator{\PGL}{PGL}

\DeclareMathOperator{\Sym}{Sym}
\newcommand{\Orbit}{\mathcal{O}}
\DeclareMathOperator{\Serre}{\mathsf{e}}
\DeclareMathOperator{\VERT}{\mathcal{V}}
\newcommand{\FF}{\mathsf{F}}
\DeclareMathOperator{\Leg}{\mathcal{L}}
\DeclareMathOperator{\Flag}{\mathcal{F}}
\DeclareMathOperator{\Edge}{\mathcal{E}}

\renewcommand{\L}{\mathsf{L}}
\newcommand{\MH}{\textsf{M}}
\newcommand{\half}{\tfrac{1}{2}}

\newcommand{\Var}{{\mathsf{V}}}
\newcommand{\XX}{{\mathcal{X}}}
\newcommand{\YY}{{\mathcal{Y}}}
\newcommand{\ZZ}{{\mathcal{Z}}}
\newcommand{\Ahat}{\widehat{A}}

\DeclareMathOperator{\Exp}{Exp}
\DeclareMathOperator{\Log}{Log}

\newcommand{\Cat}{\mathcal{C}}
\newcommand{\MM}{\mu}
\newcommand{\mm}{\nu}
\newcommand{\MMbar}{\bar{\mu}}
\newcommand{\mmbar}{\bar{\nu}}
\DeclareMathOperator{\rk}{rk}
\DeclareMathOperator{\inv}{inv}
\newcommand{\CL}{\mathcal{L}}

\begin{document}

\title{The Betti numbers of $\Mbar_{0,n}(r,d)$}

\author{Ezra Getzler}

\address{Department of Mathematics, Northwestern University}

\email{getzler@math.northwestern.edu}

\author{Rahul Pandharipande}

\address{Department of Mathematics, Princeton University}

\email{rahulp@math.princeton.edu}

\date{February 24, 2005}


\begin{abstract}
  We calculate the Betti numbers of the coarse moduli space of stable
  maps of genus 0 to projective space, using a generalization of the
  Legendre transform.
\end{abstract}

\maketitle

Let $\CM_{0,n}(r,d)$ be the moduli space of degree $d$ maps from
$n$-pointed, nonsingular, rational curves to $\P^r$ over $\C$.
Kontsevich has introduced a compactification,
\begin{equation*}
  \CM_{0,n}(r,d)\subset \Mbar_{0,n}(r,d) ,
\end{equation*}
by stable maps \cite{kont}. The moduli space $\Mbar_{0,n}(r,d)$ is a
nonsingular Deligne-Mumford stack with a projective coarse model.
Foundational issues are treated in \cite{BM, rahul}.  The main result
of our paper is a calculation of the Betti numbers of the coarse
moduli space $\Mbar_{0,n}(r,d)$.

Let $K(\Var,S_n)$ be the Grothendieck group of quasi-projective
varieties over $\C$ with action of the symmetric group $S_n$. Denote
the class of a variety $X$ in $K(\Var,S_n)$ by $[X]$. If $Y$ is an
$S_n$-invariant subspace of $X$, then
\begin{equation*}
  [X] = [X\backslash Y] + [Y] .
\end{equation*}

Let $\SS$ be the groupoid
\begin{equation*}
  \SS = \bigsqcup_{n=0}^\infty S_n .
\end{equation*}
We call a sequence of quasi-projective varieties $\XX=(\XX(n)\mid
n\ge0)$ over $\C$ with actions of $S_n$ an $\SS$-space. The
Grothendieck group of $\SS$-spaces decomposes as a product
\begin{equation*}
  K(\Var,\SS) = \prod_{n=0}^\infty K(\Var,S_n) .
\end{equation*}

Our main result is a combinatoric formula which relates the classes
\begin{equation*}
  \MM = \sum_{n=0}^\infty \sum_{d=0}^\infty q^d [\CM_{0,n}(r,d)]
  \quad \text{and} \quad
  \MMbar = \sum_{n=0}^\infty \sum_{d=0}^\infty q^d [\Mbar_{0,n}(r,d)]
\end{equation*}
of $K(\Var,\SS)\[q\]$. The formula is of a type which is now standard
in the theory of enumeration of unrooted trees, and is closely related
to the formalism of the Legendre transform for symmetric functions
introduced in \cite{modular}.

Let $K(\MH,S_n)$ be the Grothendieck group of finite-dimensional
representations of $S_n$ in the abelian category $\MH$ of mixed Hodge
structures. The (equivariant) Serre characteristic is a homomorphism
\begin{equation*}
  \Serre : K(\Var,S_n) \to K(\MH,S_n) ,
\end{equation*}
defined by the Euler characteristic with compact supports of $X$ in
$K(\MH,S_n)$:
\begin{equation*}
  \Serre(X) = \sum_i (-1)^i [H^i_c(X,\Q)] .
\end{equation*}
The Serre characteristic is a refined version of what has been
referred to as the $E$-polynomial \cite{DK} and the virtual Poincar\'e
polynomial \cite{FM}; we call it the Serre characteristic since its
existence was probably first conjectured by Serre.

If $X$ is nonsingular and projective, the weight $k$ summand of
$\Serre(X)$ is the cohomology group $H^k(X,\Q)$. The Serre
characteristic therefore generalizes the Poincar\'e polynomial of
nonsingular projective varieties.

Let $s_n\in K(\MH,S_n)$ be the class of the Hodge structure $\Q(0)$
with trivial $S_n$-action. Since $\MH$ is an abelian category over
$\Q$, we have a natural isomorphism,
\begin{equation*}
  K(\MH,\SS) \cong K(\MH)\[s_1,s_2,\dots\],
\end{equation*}
determined by induction of representations. We will give an explicit
formula for
\begin{equation*}
  \Serre(\MM)\in K(\MH)\[q,s_1,s_2,\dots\] .  
\end{equation*}
After applying a transformation analogous to the Legendre transform,
we will obtain from this our formula for
\begin{equation*}
  \Serre(\MMbar)\in K(\MH)\[q,s_1,s_2,\dots\] .
\end{equation*}

Define a morphism $\rk$ from $K(\MH,\SS)$ to $K(\MH)\[x\]$, which takes
$s_n$ to $x^n/n!$. Applying this morphism to $\Serre(\MMbar)$, we
obtain a generating function for the Serre characteristics of the
projective varieties $\Mbar_{0,n}(r,d)$:
\begin{equation*}
  \rk(\Serre(\MMbar)) = \sum_{n=0}^\infty \frac{x^n}{n!}
  \sum_{d=0}^\infty q^d \sum_{k=0}^\infty (-1)^k
  [H^k(\Mbar_{0,n}(r,d),\Q)] .
\end{equation*}
It turns out $[H^k(\Mbar_{0,n}(r,d),\Q)]$ vanishes for $k$ odd, and is
a multiple of $\L^{k/2}$ for $k$ even, where $\L=[\P]-1$, so it is
easy to determine the Betti numbers, and even the Hodge numbers, of
$H^k(\Mbar_{0,n}(r,d),\Q)$ from this result.

In the case where $r=0$, the moduli space of maps $\Mbar_{0,n}(0,d)$
is empty if $d>0$, and equals the moduli space of pointed curves
$\Mbar_{0,n}$ when $d=0$. The calculation of the Betti numbers of
$\Mbar_{0,n}$ \cite{gravity, Manin:Texel}, illustrates the use of the
Legendre transform in this problem.
 
The open stratum $\CM_{0,n}\subset\Mbar_{0,n}$ is the Zariski locally
trivial quotient of the configuration space $\FF(\P,n)$ of the
projective line $\P$ by its automorphism group $\Aut(\P)$, hence
\begin{equation*}
  \Serre(\CM_{0,n}) = \frac{\Serre(\FF(\P,n))}{\Serre(\Aut(\P))} =
  \frac{\Serre(\FF(\P,n))}{\L(\L^2-1)} .
\end{equation*}
%\begin{equation*}
%  \Serre(\CM_{0,n}) = \frac{\Serre(\FF(\P,n))}{\Serre(\Aut(\P))}.
%\end{equation*}
Introduce the power series
\begin{equation*}
  \mm(x) = \sum_{n=3}^\infty \frac{x^n}{n!} \rk(\Serre(\CM_{0,n}))
  \quad \text{and} \quad
  \mmbar(x) = \sum_{n=3}^\infty \frac{x^n}{n!}
  \rk(\Serre(\Mbar_{0,n})) .
\end{equation*}
We have a stratification
\begin{equation*}
  \Mbar_{0,n} = \sum_{T\in\Gamma_{0,n}} \prod_{v\in\VERT(T)} \CM_{0,n(v)} ,
\end{equation*}
where $\Gamma_{0,n}$ is the set of isomorphism classes of trees $T$,
all vertices $v\in\VERT(T)$ of which have valence $n(v)>2$, with $n$
labelled leaves. Then $\half x^2-\mm$ and $\half x^2+\mmbar$ are
Legendre transforms of each other: if $D$ is the operation of
differentiation with respect to $x$, then
\begin{equation*}
  \mmbar = \mm\circ(x+D\mmbar) - \half(D\mmbar)^2 .
\end{equation*}
Applying $D$ to both sides of this equation, we see that
\begin{equation*}
  D\mmbar=D\mm\circ(x+D\mmbar) ,  
\end{equation*}
in other words, that the series $x-D\mm$ and $x+D\MMbar$ are inverse
to each other:
\begin{equation*}
  (x-D\mm)\circ(x+D\mmbar) = x .
\end{equation*}
This allows the recursive calculation of $D\mmbar$, and hence of
$\mmbar$.

The moduli space of maps to $\P^r$ has a similar stratification
indexed by trees. However, for $r>0$, the trees may have non-trivial
automorphisms. The strata are naturally described as \emph{quotients}
of fibered products of vertex moduli spaces by the group of tree
automorphisms.  Our formulas for $\P^r$ require knowledge of the
equivariant Serre characteristics of the spaces $\CM_{0,n}(r,d)$; it
is for this reason that we work with the Grothendieck groups
$K(\Var,\SS)$ and $K(\MH,\SS)$, and not just $K(\Var)$ and $K(\MH)$.

Our methods also calculate the Betti numbers of the moduli spaces
$\Mbar_{0,n}(r,d)/S_n$ of stable maps with unlabelled marked points:
if
\begin{equation*}
  \inv:K(\MH,\SS)\to K(\MH)\[x\]
\end{equation*}
is defined by mapping $s_n$ to $x^n$, then we have
\begin{equation*}
  \inv(\Serre(\MMbar)) = \sum_{n=0}^\infty {x^n} \sum_{d=0}^\infty
  q^d \sum_{k=0}^\infty (-1)^k [H^k(\Mbar_{0,n}(r,d)/S_n,\Q)] .
\end{equation*}

\subsection*{Related work}

The moduli spaces $\Mbar_{0,n}(r,d)$ are rational \cite{bp} with
algebraic cohomology. In fact, the cohomology groups are generated by
tautological classes \cite{dr1, dr2}.  The Picard group of
$\Mbar_{0,n}(r,d)$ is determined in \cite{rint}.  Betti number
calculations in degree 2 can be found in \cite{cox}.  The
cohomological ring structure has been studied in degree 2 (with
partial results in degree 3) in \cite{bo}. A more general approach to
the ring structure has been recently proposed in \cite{mus}.

\subsection*{Acknowledgements}

We are grateful to B. Totaro for suggesting the argument for Theorem
\ref{almost-free}. We thank B. Sudakov for a helpful conversation
about trees.

While writing the paper, we were partially supported by the NSF and
the Packard Foundation. The paper was started in 1995 at the Miracle
of Science in Cambridge, MA.  The first draft was written while E.G.
was a member of the Max-Planck Institute in Bonn and R.P. was a member
of the Mittag-Leffler Institute in Stockholm. The paper was finished
in Princeton while E.G. was visiting the Institute for Advanced Study
and Princeton University.

\section{Composition algebras}

Consider the algebra of polynomials $A=R[x]$ in one variable over a
commutative ring $R$. This algebra is filtered by subalgebras
\begin{equation*}
  F_nA = x^nA .
\end{equation*}
There is a composition operation $\circ:A\o F_1A\to A$, characterized
by the following axioms:
\begin{enumerate}
\item for fixed $b\in F_1A$, the map $a\mapsto a\circ b$ is an element
  of the endomorphism algebra $\text{End}_R(A)$, and
  \begin{enumerate}
  \item $a\circ b = a$ if $a\in R$,
  \item $(a_1+a_2)\circ b = a_1\circ b+a_2\circ b$,
  \item $(a_1a_2)\circ b = (a_1\circ b)(a_2\circ b)$;
  \end{enumerate}
\item $(a\circ b_1)\circ b_2=a\circ(b_1\circ b_2)$;
\item $x\circ b = b$ and $a \circ x = a$.
\end{enumerate}
There is a derivation
\begin{equation*}
  Da = \frac{da}{dx}
\end{equation*}
of the algebra $A$ over $R$, such that $Dx=1$, and
\begin{equation*}
  D(a\circ b)=(Da\circ b)(Db) .
\end{equation*}

Suppose that $R$ has characteristic $0$. Define element $s_n$ of $A$
by
\begin{equation*}
  s_n = \frac{x^n}{n!} .
\end{equation*}
The relations
\begin{equation*}
  s_n(b_1+b_2) = \sum_{k=0}^n (s_k\circ b_1)(s_{n-k}\circ b_2)
\end{equation*}
and $Ds_n=s_{n-1}$ hold.

The goal of this section is to axiomatize the composition operations
of the above type, and derive some elementary consequences.
\begin{definition}
  Let $A$ be a filtered commutative algebra over a commutative ring
  $R$:
  \begin{equation*}
    A = F_0A \supset F_1A \supset F_2A \supset \dots .
  \end{equation*}
  A composition operation on $A$ is an operation $\circ:A\o F_1A\to A$
  such that
  \begin{equation*}
    F_kA\circ F_\ell A\subset F_{k\ell}A ,
  \end{equation*}
  together with a sequence of elements $s_n\in F_nA$ for $n\ge0$, and
  a derivation
  \begin{equation*}
    D:F_nA\to F_{n-1}A ,    
  \end{equation*}
  satisfying the following axioms:
  \begin{enumerate}
  \item for fixed $b\in F_1A$, the map $a\mapsto a\circ b$ is an
    element of the endomorphism algebra $\text{End}_R(A)$;
  \item $(a\circ b_1)\circ b_2=a\circ(b_1\circ b_2)$;
  \item $D(a\circ b)=(Da\circ b)(Db)$;
  \item $s_0=1$, $s_1\circ b=b$, $a\circ s_1=a$, $Ds_{n}=s_{n-1}$, and
    \begin{equation*}
      s_n(b_1+b_2) = \sum_{k=0}^n (s_k\circ b_1)(s_{n-k}\circ b_2) .
    \end{equation*}
  \end{enumerate}
\end{definition}

If $A$ is an algebra with composition operation, then so is the
completion $\Ahat$ with respect to the filtration $F_\bull A$.

\section{The composition operation on the Grothendieck group of
  $\SS$-varieties}

The Grothendieck group $K(\Var)$ of quasi-projective varieties over
$\C$ is generated by elements $[X]$, where $X$ is a quasi-projective
variety over $\C$, subject to the following relations:
\begin{enumerate}
\item $[X]=[Y]$ if $X$ and $Y$ are isomorphic as varieties;
\item $[X] = [X\backslash Y]+[Y]$ when $Y$ is a subvariety of $X$.
\end{enumerate}
The Cartesian product of varieties induces a commutative product on
$K(\Var)$, with identity $1=[\Spec(\C)]$.

Let $\Gamma$ be a finite group. The Grothendieck group
$K(\Var,\Gamma)$ is generated by elements $[X]$, where $X$ is a
quasi-projective variety over $\C$ with action of $\Gamma$, subject to
the following relations:
\begin{enumerate}
\item $[X]=[Y]$ if $X$ and $Y$ are isomorphic as varieties with action
  of $\Gamma$;
\item $[X] = [X\backslash Y]+[Y]$ when $Y$ is a $\Gamma$-invariant
  subvariety of $X$.
\end{enumerate}
The Cartesian product of varieties makes $K(\Var,\Gamma)$ into a
$K(\Var)$-module.

\begin{proposition}
  Let $\XX$ be a locally trivial $\Gamma$-equivariant fibration over a
  variety $B$ with trivial action of $\Gamma$, and with fibre $\YY$.
  Then the elements $[\XX]$ and $[B]\,[\YY]$ of $K(\Var,\Gamma)$ are
  equal.
\end{proposition}
\begin{proof}
  Take a Zariski open subset $U\subset B$ such that
  $\pi^{-1}(U)\cong\YY\times U$. We then have
  \begin{equation*}
    [\XX] = [\pi^{-1}(U)] + [\pi^{-1}(B\setminus U)] = [\YY] \, [U] +
    [\pi^{-1}(B\setminus U)] .
  \end{equation*}
  The proposition follows by Noetherian induction on $B$.
\end{proof}

There is an analogue of the Grothendieck group of a finite group
$\Gamma$, where $\Gamma$ is replaced by the groupoid $\SS$. The
resulting Grothendieck group decomposes as a product
\begin{equation*}
  K(\Var,\SS) = \prod_{n=0}^\infty K(\Var,S_n) .
\end{equation*}
The $\boxtimes$-product on $\SS$-spaces,
\begin{equation*}
  \bigl( \XX \boxtimes \YY \bigr)(n) = \bigsqcup_{m=0}^n
  \Ind_{S_m\times S_{n-m}}^{S_n} \bigl( \XX(m) \times \YY(n-m) \bigr) ,
\end{equation*}
induces a commutative product on $K(\Var,\SS)$. In particular,
$K(\Var,\SS)$ is an algebra over the commutative ring
\begin{equation*}
  K(\Var) = K(\Var,S_0) \subset K(\Var,\SS) .
\end{equation*}

Let $F_nK(\Var,\SS)$ be the subspace of $K(\Var,\SS)$ spanned by
classes $[\XX]$ such that $\XX(m)$ is empty for $m<n$. These subspaces
define a filtration of the algebra $K(\Var,\SS)$.

Define the composition of $\SS$-spaces $\XX$ and $\YY$ by
\begin{equation*}
  \bigl( \XX\circ\YY \bigr) (n) = \bigsqcup_{i=0}^\infty
 \bigl( \XX(i) \times \YY^{\boxtimes i}(n) \bigr) \big/ S_i .
\end{equation*}
This operation is associative:
\begin{equation*}
  (\XX\circ\YY)\circ\ZZ \cong \XX\circ(\YY\circ\ZZ) .
\end{equation*}
There is a unique associative operation $\circ$ on $K(\Var,\SS)\times
F_1K(\Var,\SS)$ such that
\begin{equation*}
  [\XX]\circ[\YY] = [\XX\circ\YY].
\end{equation*}
We leave the straightforward verification to the reader.

There is a functor $\XX\mapsto \delta\XX$ from $\SS$-spaces to
$\SS$-spaces, defined by
\begin{equation*}
  \delta\XX(n) = \Res^{S_{n+1}}_{S_n} \XX(n+1). 
\end{equation*}
It is easily seen that
\begin{equation*}
  \delta(\XX\boxtimes\YY) \cong \delta\XX\boxtimes\YY \sqcup
  \XX\boxtimes \delta\YY .
\end{equation*}
This functor induces a derivation $D$ of the algebra $K(\Var,\SS)$
over $K(\Var)$, by the formula $\delta[\XX]=[\delta\XX]$, which
satisfies
\begin{equation*}
  D:F_nK(\Var,\SS) \to F_{n-1}K(\Var,\SS) .
\end{equation*}
Since
\begin{equation*}
  \delta\bigl( \XX^{\boxtimes i} \bigr) \cong \Ind_{S_{i-1}}^{S_i}
  \XX^{\boxtimes i-1} \boxtimes \delta\XX ,
\end{equation*}
we also see that
\begin{align*}
  \delta(\XX\circ\YY) &= \bigsqcup_{i=0}^\infty \bigl( \XX(i) \times
  \delta\bigl( \YY^{\boxtimes i} \bigr) \bigr) \big/ S_i \\
  &= \bigsqcup_{i=1}^\infty \bigl( \XX(i) \times \Ind_{S_{i-1}}^{S_i}
  \bigl( \YY^{\boxtimes i-1} \boxtimes \delta\YY \bigr) \bigr) \big/ S_i \\
  &= \bigsqcup_{i=0}^\infty \bigl( \delta\XX(i) \times \bigl(
  \YY^{\boxtimes i} \boxtimes \delta\YY \bigr) \bigr) \big/ S_i ,
\end{align*}
and hence that $D$ satisfies the equation
\begin{equation*}
  D(a\circ b) = (Da\circ b) \, Db .
\end{equation*}

Let $s_n$ be the class in $K(\Var,S_n)\subset K(\Var,\SS)$ associated
to the $\SS$-space $\Spec(\C)$ with trivial action of $S_n$. We have
$s_0=1$,
\begin{equation*}
  s_1\circ[\XX] = [\XX]\circ s_1 = [\XX] ,
\end{equation*}
and $Ds_n=s_{n-1}$. In particular, $Ds_1=1$.

Assembling the above constructions on $K(\Var,\SS)$, we obtain the
following theorem.
\begin{theorem} \label{composition:motive}
  The algebra $K(\Var,\SS)$ is a complete algebra with composition
  operation over $K(\Var)$.
\end{theorem}

\section{Some calculations in $K(\Var,\SS)$}

In this section, we present some formulas in $K(\Var,\SS)$ which will
be needed in the formulation of our results concering the Betti
numbers of $\Mbar_{0,n}(r,d)$.

Let $\Sigma$ be the universal algebra with composition operation, with
generators $s_n\in F_n\Sigma$. Any complete algebra with composition
operation carries an analogue of the exponential function $e^x-1$,
\begin{equation*}
  \Exp(a) = \sum_{n=1}^\infty s_n\circ a : F_1A\to F_1A ,
\end{equation*}
which satisfies the product formula
\begin{equation} \label{Cartan}
  \Exp(a+b) = \Exp(a)\Exp(b) + \Exp(a) + \Exp(b) .
\end{equation}
For example, let  $X$ be is a quasi-projective variety over $\C$,
and let $[X] \in K(\Var,S_0)$. Then,
\begin{equation*}
  \Exp([X]s_1)(n) =
  \begin{cases}
    0 , & n=0 , \\
    [X^n] , & n>0 ,
  \end{cases}
\end{equation*}
where $S_n$ acts on the $S_n$-space $X^n$ by permuting the factors.

If $A$ is a complete algebra with composition operation, the operation
$\Exp$ has an inverse $\Log$.
\begin{proposition}
  There is a sequence of elements $\ell_n\in F_n\Sigma$, such that the
  operation
  \begin{equation*}
    \Log(a) = \sum_{n=1}^\infty \ell_n\circ a : F_1A \to F_1A
  \end{equation*}
  is the inverse of $\Exp$. For example, $\ell_1=s_1$, $\ell_2=-s_2$,
  $\ell_3=-s_3+s_2s_1$ and
  \begin{equation*}
    \ell_4 = - s_4 + s_3s_1 + s_2\circ s_2 - s_2s_1^2 .
  \end{equation*}
\end{proposition}
\begin{proof}
  The equation $\Exp(\Log(a))=a$ is equivalent to
  \begin{equation*}
    \sum_{n,k=1}^\infty s_n \circ \ell_k  = s_1 .
  \end{equation*}
  This implies that $\ell_1=s_1$, while for $n>1$,
  \begin{equation*}
    \sum_{k_1+2k_2+\dots nk_n=n} (s_{k_1}\circ\ell_1) \dots
    (s_{k_n}\circ \ell_n) = 0 .
  \end{equation*}
  This equation determines $\ell_n$ in terms of $\ell_k$, $k<n$.
\end{proof}

The configuration space $\FF(X)$ of a variety $X$ is the $\SS$-space
consisting of embeddings of the discrete variety $\{1,\dots,n\}$ with
$n$ points into $X$:
\begin{equation*}
  \FF(X,n) = \{ (z_1,\dots,z_n) \in X^n \mid \text{$z_i\ne z_j$ for
  $i\ne j$} \} .
\end{equation*}
Then $\FF(X,n)$ is a quasi-projective variety on which $S_n$ acts by
permutation of the points. Define the $\SS$-space $\FF(X)$ by 
\begin{equation*}
  \FF(X)(n) = \FF(X,n) .
\end{equation*}
The following formula for $[\FF(X)]$ is a generalization of a formula
of \cite{I}, and is proved by the same method.
\begin{theorem}
  $[\FF(X)] = 1+\Exp\bigl( [X] \Log(s_1) \bigr)$  
\end{theorem}
\begin{proof}
  The space $X^n$ has an $S_n$-equivariant decomposition into locally
  closed subvarieties
  \begin{align*}
    X^n_i &= \{ (x_1,\dots,x_n) \in X^n \mid \text{the set
      $\{x_1,\dots,x_n\}$ has $i$ distinct points} \} \\
    &= \bigsqcup_{\substack{n_1+\dots+n_i=n \\ n_1,\dots,n_i>0}}
    \bigl( \FF(X,i) \times \Ind^{S_n}_{S_{n_1}\times\dots\times
      S_{n_i}} \Spec(\C) \bigr) \big/ S_i .
  \end{align*}
  Taking the union over $1\le i\le n$, we see that
  \begin{equation*}
    \Exp([X]s_1) = [\FF(X)] \circ \Exp(s_1) .
  \end{equation*}
  It follows that
  \begin{align*}
    \Exp([X]s_1) \circ \Log(s_1) &= [\FF(X)] \circ \Exp(s_1) \circ
    \Log(s_1) \\ &= [\FF(X)] \circ s_1 = [\FF(X)] .
  \end{align*}
  On the other hand, we have $\Exp([X]s_1) \circ \Log(s_1) =
  \Exp\bigl( [X] \Log(s_1) \bigr)$.
\end{proof}

A graded $\SS$-space is a sequence $(\XX(n,d)\mid d\ge0)$ of
$\SS$-spaces. The associated Grothendieck group is naturally
isomorphic to $K(\Var,\SS)\[q\]$, by the identification
\begin{equation*}
  [\XX] = \sum_{d=0}^\infty q^d [\XX(-,d)] .
\end{equation*}
Filter $K(\Var,\SS)$ by subspaces
\begin{equation*}
  F_nK(\Var,\SS)\[q\] = \bigcup_{i=0}^n q^i \bigl(
  F_{n-i}K(\Var,\SS) \bigr)\[q\] .
\end{equation*}

Define the composition of graded $\SS$-spaces $\XX$ and $\YY$ by
\begin{equation*}
  \bigl( \XX\circ\YY \bigr) (n,d) = \bigsqcup_{i=0}^\infty
  \bigsqcup_{e=0}^d \bigl( \XX(i,d-e) \times \YY^{\boxtimes i}(n,e)
  \bigr) \big/ S_i .
\end{equation*}
This induces an associative operation $\circ$ on
$K(\Var,\SS)\[q\]\times F_1K(\Var,\SS)\[q\]$, which makes
$K(\Var,\SS)\[q\]$ into a complete algebra with composition
operation over $R=K(\Var)\[q\]$.

Let $\Map(r,d)=\Map_d(\P,\P^r)$ be the variety of algebraic maps of
degree $d$ from the projective line $\P$ to $\P^r$. The following is
the main result of this section.
\begin{proposition}
  \begin{equation*}
    [\Map(r,d)] =
    \begin{cases} 
      [\P^r] , & d=0 , \\
      [\A]^{(d-1)(r+1)+1} ([\A]^r-1) [\P^r] , & d>0 .
    \end{cases}
  \end{equation*}
\end{proposition}
\begin{proof}
  Let $V$ be a two-dimensional vector space over $\C$. An element of
  the space $\Map_d(\P(V),\P^r)$ of algebraic maps of degree $d$ from
  the projective line $\P(V)$ to $\P^r$ corresponds to a sequence
  $$(f_0,\dots,f_r)\in\Sym^d(V^*)\o\C^{r+1}$$ 
of homogeneous polynomials
  of degree $d$ with no common roots: there is an open
  embedding
  \begin{equation*}
    \Map_d(\P(V),\P^r) \hookrightarrow \P(\Sym^d(V^*)\o\C^{r+1}) ,
  \end{equation*}
  whose complement is the resultant consisting of sequences of
  polynomials with a common root.
  
  Stratify $\P(\Sym^d(V^*)\o\C^{r+1})$ by the number of common roots
  of the $r+1$ polynomials $(f_0,\dots,f_r)$ to obtain a decomposition
  \begin{equation*}
    \P(\Sym^d(V^*)\o\C^{r+1}) = \bigsqcup_{e=0}^d \Map_e(\P(V),\P^r)
    \times \P(\Sym^{d-e}(V^*)) .
  \end{equation*}
  Taking the class in $K(\Var)$ of both sides, we see that
  \begin{align*}
    \sum_{d=0}^\infty q^d [\Map(r,d)] &=
    \frac{\displaystyle\sum_{d=0}^\infty q^d
      [\P(\Sym^d(V^*)\o\C^{r+1})]} {\displaystyle \sum_{d=0}^\infty
      q^d [\P(\Sym^d(V^*))]} \\
    &= \frac{\displaystyle\sum_{d=0}^\infty q^d [\P^{(d+1)(r+1)-1}]}
    {\displaystyle \sum_{d=0}^\infty q^d [\P^d]} .
  \end{align*}
  We have
  \begin{align*}
    \sum_{d=0}^\infty q^d [\P^{(d+1)(r+1)-1}] &= \frac{1}{[\A]-1}
    \sum_{d=0}^\infty q^d ([\A]^{(d+1)(r+1)}-1) \\
    &= \frac{1}{[\A]-1} \left( \frac{[\A]^{r+1}}{1-q[\A]^{r+1}} -
      \frac{1}{1-q} \right) \\
    &= \frac{[\P^r]}{(1-q)(1-q[\A]^{r+1})} ,
  \end{align*}
  while
  \begin{equation*}
    \sum_{d=0}^\infty q^d [\P^d] = \frac{1}{(1-q)(1-q[\A])} .
  \end{equation*}
  In this way, we see that
  \begin{equation*}
    \sum_{d=0}^\infty q^d [\Map(r,d)] = \frac{1-q[\A]}{1-q[\A]^{r+1}}
    \, [\P^r] .
  \end{equation*}
  The proposition follows on expanding the power series on the
  right-hand side.
\end{proof}

\section{The moduli space $\Mbar_{0,n}(r,d)$ of stable maps}

Our calculation of $[\Mbar_{0,n}(r,d)]$ depends on the stratification
of $\Mbar_{0,n}(r,d)$ by strata indexed by stable trees, which we now
recall (see \cite{manin} for further details).

\begin{definition}
  A graph $T$ is given by the data $(\Flag(T),\sigma,\tau)$, where
  $\Flag(T)$ is the set of flags of $T$, $\sigma$ is an involution of
  $\Flag(T)$, and $\tau$ is an equivalence relation on $\Flag(T)$.
\end{definition}

The geometric realization $|T|$ of a graph $T$ is the one-dimensional
cellular complex, constructed as follows.
\begin{itemize}
\item The leaves of $T$ form the boundary of $|T|$; they are elements
  of the set $\Leg(T)$ of fixed points of the involution $\sigma$.
\item The vertices of $T$ are the $0$-cells of $|T|$; they are
  elements of the set $\VERT(T)$ of equivalence classes of $\Flag(T)$
  with respect to $\tau$.
\item The edges of $T$ are the $1$-cells of $|T|$ joining two
  vertices; they are elements of the set $\Edge(T)$ of orbits of
  $\sigma$ with two elements.
\item The remaining $1$-cells of $|T|$ join a leaf to the vertex in
  whose $\tau$-equivalence class it lies, and are in bijective
  correspondence with the set $\Leg(T)$.
\end{itemize}
The subset of $\Flag(T)$ which meet a vertex $v\in\VERT(T)$ is denoted
$\Leg(v)$, and its cardinality $n(v)$ is called the valence of $v$.

A tree $T$ is a graph whose geometric realization is simply connected,
that is, $|T|$ is connected and
\begin{equation*}
   |\Flag(T)| = |\VERT(T)| + |\Edge(T)|+ |\Leg(T)| - 1 .
\end{equation*}

An $n$-tree $T$, where $n>0$, is a tree $T$ together with a bijection
between its leaves and the set $\{1,\dots,n\}$. An $(n,d)$-tree $T$,
where $d\ge0$, is an $n$-tree $T$, together with a function $v\mapsto
d(v)$ from the vertices of $T$ to $\N$, such that
\begin{equation*}
  \sum_{v\in\VERT(T)} d(v) = d .
\end{equation*}
An automorphism of an $(n,d)$-tree $T$ is an automorphism of the
underlying tree, that is, an automorphism of $\Flag(T)$ compatible
with the involution $\sigma$ and partition $\tau$, which fixes the
labelling of its leaves and preserves the function $d(v)$. The group of
automorphisms is denoted $\Aut(T)$.

\begin{definition}
  An $(n,d)$-tree is stable if for each vertex $v\in\VERT(T)$ either
  $d(v)>0$ or $n(v)>2$.
\end{definition}
Denote the set of isomorphism classes of stable $(n,d)$-trees by
$\Gamma_{0,n}(d)$.

\begin{proposition}
  The set $\Gamma_{0,n}(d)$ of stable $(n,d)$-trees is finite for each
  $n$ and $d$.
\end{proposition}
\begin{proof}
  Since the geometric realization of a tree $T$ is simply-connected
  \begin{equation*}
    \sum_{v\in\VERT(T)} (n(v)-2) = n-2 .
  \end{equation*}
  The number of vertices for which $n(v)\le2$ is bounded by $d$, while
  the number of vertices for which $n(v)>2$ is bounded by $n-2$. Hence
  the number of flags is bounded by $3n+2d-4$. The number of trees
  with a fixed number of flags is finite, proving the proposition.
\end{proof}

A stable map $(f:\Sigma\to\P^r,z_1,\dots,z_n)$ of genus $0$ consists
of the following data.
\begin{itemize}
\item A complete, connected  curve $\Sigma$ of arithmetic genus $0$ whose only
  singularities are double points.
\item Distinct marked points $(z_1,\dots,z_n)$ lying in the smooth
  locus of $\Sigma$.
\item An algebraic map $f:\Sigma\to\P^r$ with a finite number of
  automorphisms fixing the marked points.
\end{itemize}

Let $H\in H_2(\P^r)$ be the class of a line. The degree of the stable
map is the natural number $d$, where $f_*[\Sigma]=dH\in H_2(\P^r)$.

The dual graph $T$ of a stable map of genus $0$ is defined as follows.
\begin{itemize}
\item $T$ has one vertex for each irreducible component of $\Sigma$.
\item Each double point of $\Sigma$ corresponds to an edge of the dual
  tree, which joins the corresponding vertices.
\item Each marked point of $\Sigma$ corresponds to a leaf of the dual
  tree, which is attached to the vertex corresponding to the unique
  irreducible component on which the marked point lies. 
\item Given a vertex $v\in\VERT(T)$, denote by $\Sigma(v)$ the
  corresponding irreducible component of $\Sigma$, and associate to
  $v$ the degree $d(v)\in\N$ of $f$ on the irreducible component
  $\Sigma(v)$.
\end{itemize}
Since $\Sigma$ has arithmetic genus $0$, we see that $T$ is an
$(n,d)$-tree.

The following proposition is an immediate consequence of the definitions of
stable trees and stable maps.
\begin{proposition}
  A map $(f:\Sigma\to\P^r,z_1,\dots,z_n)$ of genus $0$ with $n$ marked
  points and degree $d$ is stable if and only if its dual graph is a
  stable $(n,d)$-tree.
\end{proposition}

The moduli stack of stable maps of genus $0$ with $n$ marked points
and degree $d$ to $\P^r$ is a smooth, complete Deligne-Mumford stack
The associated coarse moduli
space, $\Mbar_{0,n}(r,d)$, is a projective variety whose only
singularities are quotient singularities by finite groups. We refer the
reader to \cite{BM}, \cite{rahul} for a treatment of foundational issues.

Given $T\in\Gamma_{0,n}(d)$, let $\CM(T)$ be the locally closed
subvariety of $\Mbar_{0,n}(r,d)$ consisting of the moduli of all
stable maps with dual graph $T$. 
The codimension of $\CM(T)$ is $|\Edge(T)|$. The
collection of subvarieties $\CM(T)$ is a stratification of
$\Mbar_{0,n}(r,d)$,
$$[\Mbar_{0,n}(r,d)] = \sum_{T\in \Gamma_{0,n}(d)} [\CM(T)]\ \  \in
K(\Var,S_n).
 $$
The right side of the above formula represents a slight abuse of terminology
since only the $S_n$-invariant unions of strata define
elements of $K(\Var,S_n)$.

There is a unique open stratum $\CM_{0,n}(r,d)$ of $\Mbar_{0,n}(r,d)$,
corresponding to the unique stable tree with a single vertex of
valence $n$ (and no edges) and
 equal to the quotient
\begin{equation*}
  \CM_{0,n}(r,d) = \bigl( \FF(\P,n)\times\Map_d(\P,\P^r) \bigr) \big/
  \Aut(\P) .
\end{equation*}

If $\CL$ is an arbitrary finite set, let $\CM_{0,\CL}(r,d)$ denote the
moduli space of maps of degree $d$ from $\P$ to $\P^r$ together with
an embedding of the finite set $\CL$ into $\P$. The moduli space
$\CM_{0,\CL}(r,d)$ is isomorphic to $\CM_{0,|\CL|}(r,d)$, but not in
any natural way. A canonical evaluation map
\begin{equation*}
  \CM_{0,\CL}(r,d) \rightarrow (\P^r)^\CL   
\end{equation*}
is given by the composition of the embedding $\CL\hookrightarrow\P$
with the map $f:\P\to\P^r$. Since the set of flags $\Flag(T)$ of a
tree $T$ is the disjoint union of the sets $\Leg(v)$, we obtain a map
\begin{equation*}
  \prod_{v\in\VERT(T)} \CM_{0,\Leg(v)}(r,d(v)) \to (\P^r)^{\Flag(T)} .
\end{equation*}
There is also a natural map from
$\bigl(\P^r\bigr)^{\Edge(T)\sqcup\Leg(T)}$ to
$\bigl(\P^r\bigr)^{\Flag(T)}$, which may be thought of as the
inclusion of the fixed-point set under the action of the involution
$\sigma$ on $\bigl(\P^r\bigr)^{\Flag(T)}$. The stratum $\CM(T)$ is
naturally isomorphic to the quotient of the fibred product
\begin{equation} \label{CM}
  \begin{CD}
    \CM_ \square(T) @>>> \prod_{v\in\VERT(T)} \CM_{0,\Leg(v)}(r,d(v)) \\
    @VVV @VVV \\
    \bigl(\P^r\bigr)^{\Edge(T)\sqcup\Leg(T)} @>>> (\P^r)^{\Flag(T)}
  \end{CD}
\end{equation}
by the finite group $\Aut(T)$,
\begin{equation*}
  \CM(T)= \CM_\square(T)/\Aut(T) .  
\end{equation*}

Let $\CM^*_{0,n}(r,d)$ and $\Mbar^*_{0,n}(r,d)$ be the fibres of the
evaluation maps
\begin{align*}
  & (f:\Sigma\to\P^r,z_1,\dots,z_{n+1}) \mapsto f(z_{n+1}) :
  \CM_{0,n+1}(r,d) \to \P^r \\
  & (f:\Sigma\to\P^r,z_1,\dots,z_{n+1}) \mapsto f(z_{n+1}) :
  \Mbar_{0,n+1}(r,d) \to \P^r.
\end{align*}
Both these fibrations are Zariski locally trivial, by elementary
considerations.

Define elements of $K(\Var,\SS)\[q\]$ by
\begin{equation*}
  [\CM(r)] = \sum_{n=0}^\infty \sum_{d=0}^\infty q^d [\CM_{0,n}(r,d)],
  \quad \quad
  [\CM^*(r)] = \sum_{n=0}^\infty \sum_{d=0}^\infty q^d [\CM^*_{0,n}(r,d)]
\end{equation*}
\begin{equation*}
  [\Mbar(r)] = \sum_{n=0}^\infty \sum_{d=0}^\infty q^d [\Mbar_{0,n}(r,d)],
  \quad \quad
  [\Mbar^*(r)] = \sum_{n=0}^\infty \sum_{d=0}^\infty q^d [\Mbar^*_{0,n}(r,d)]
\end{equation*}
By the definition of $D$ and the Zariski local triviality of the
evaluation fibrations, we have
\begin{equation*}
  D[\CM(r)] = [\P^r]\,[\CM^*(r)] \quad \text{and} \quad
  D[\Mbar(r)] = [\P^r]\,[\Mbar^*(r)] .  
\end{equation*}

The main geometric result of the paper is a relationship between the
elements $[\CM(r)]$, $[\Mbar(r)]$, and $[\Mbar^*(r)]$ in the
composition algebra $K(\Var,\SS)\[q\]$.
\begin{theorem} \label{legendre}
  $[\Mbar(r)] = [\CM(r)] \circ ( s_1 + [\Mbar^*(r)]) + [\P^r] \bigl(
  s_2 \circ [\Mbar^*(r)] - [\Mbar^*(r)]^2 \bigr)$
\end{theorem}

Applying the derivation $D$ to both sides of Theorem \ref{legendre}
and using the formula $D(s_2\circ a) = a\,Da$, we see that
\begin{equation*}
  [\P^r] \bigl( 1 + D[\Mbar^*(r)] \bigr) \bigl( [\Mbar^*(r)] -
  [\CM^*(r)] \circ ( s_1 + [\Mbar^*(r)] ) \bigr) = 0 .
\end{equation*}
Since $[\P^r]\bigl(1+D[\Mbar^*(r)]\bigr)$ is not a zero divisor in
$K(\Var,\SS)\[q\]$, we obtain the following corollary.
\begin{corollary}
  $[\Mbar^*(r)] = [\CM^*(r)] \circ \bigl( s_1 + [\Mbar^*(r)] \bigr)$
\end{corollary}

This corollary gives a recursive algorithm for calculating
$[\Mbar^*(r)]$ in terms of $[\CM^*(r)]$. Substituting the resulting
formula for $[\Mbar^*(r)]$ into Theorem \ref{legendre}, we obtain an
algorithm for calculating $[\Mbar(r)]$ from $[\CM(r)]$.

The proof of Theorem \ref{legendre} occupies the remainder of this
section. The following lemma follows easily from the definition of the
composition operation and \eqref{CM}.
\begin{lemma} \label{ccc}
  Consider the sets of pairs
  \begin{align*}
    \Gamma_{0,n}^{\VERT{}}(d) &= \{ (T,v) \mid T \in \Gamma_{0,n}(d) ,
    v \in \VERT(T) \} , \\
    \Gamma_{0,n}^{\Edge{}}(d) &= \{ (T,e) \mid T \in \Gamma_{0,n}(d) ,
    e \in \Edge(T) \} , \\
    \Gamma_{0,n}^{\Leg{}}(d) &= \{ (T,i) \mid T \in \Gamma_{0,n}(d) ,
    i \in \Leg(T) \} , \\
    \Gamma_{0,n}^{\Flag{}}(d) &= \{ (T,f) \mid T \in \Gamma_{0,n}(d) ,
    f \in \Flag(T) \} .
  \end{align*}
  Then
  \begin{align*}
    \sum_{n=0}^\infty \sum_{(T,v)\in\Gamma_{0,n}^{\VERT{}}(d)}
    [\CM_\square(T)/ \Aut(T,v)] &= [\CM(r)]\circ(s_1+[\Mbar^*(r)]) , \\
    \sum_{n=0}^\infty \sum_{(T,e)\in\Gamma_{0,n}^{\Edge{}}(d)}
    [\CM_\square(T)/ \Aut(T,e)] &= [\P^r] \, (s_2\circ[\Mbar^*(r)]) , \\
    \sum_{n=0}^\infty \sum_{(T,i)\in\Gamma_{0,n}^{\Leg{}}(d)}
    [\CM_\square(T)/ \Aut(T,i)] &=  s_1 D[\Mbar(r)] , \\
    \sum_{n=0}^\infty \sum_{(T,f)\in\Gamma_{0,n}^{\Flag{}}(d)}
    [\CM_\square(T)/ \Aut(T,f)] &= [\P^r] \, [\Mbar^*(r)]^2+s_1D[\Mbar(r)] .
  \end{align*}
\end{lemma}

\begin{lemma}\label{treel}
  Let $T\in \Gamma_{0,n}(d)$. There is a canonical injection
  \begin{equation*}
    \iota : \Flag(T) \hookrightarrow \VERT(T) \sqcup \Edge(T) \sqcup \Leg(T).
  \end{equation*}
\end{lemma}
\begin{proof}
  Define $\iota$ on the subset $\Leg(T) \subset \Flag(T)$ to be the
  identity map from this subset to the copy of $\Leg(T)$ on the right.
  Define $T^-$ by removing all legs from $T$. To complete the proof,
  we will construct a canonical injection
  \begin{equation*}
    \iota : \Flag(T^-) \hookrightarrow \VERT(T^-) \sqcup \Edge(T^-)
  \end{equation*}
  Note that there are natural identifications
  $\VERT(T^-)\cong\VERT(T)$ and $\Edge(T^-)\cong\Edge(T)$.
  
  For each extremal vertex $v$ of $T^-$ (vertex of valence $1$), with
  extremal flag $f$, let
  \begin{equation*}
    \iota(f) = v .
  \end{equation*}
  For each extremal edge $e$, made up of an extremal flag $f$ and a
  non-extremal flag $f'$, let
  \begin{equation*}
    \iota(f') = e .
  \end{equation*}
  Prune the extremal edges and vertices of $T^-$ and repeat.
\end{proof}

Since $|\Flag(T)|= |\VERT(T)|+|\Edge(T)|+|\Leg(T)|-1$, the complement
of $\Flag(T)$ under $\iota$ consists of a canonical element $\ast$ of
$\VERT(T)\sqcup \Edge(T)$.  Lemma \ref{treel} defines a canonical
bijection
\begin{equation} \label{zzz}
  \ast \sqcup \Flag(T) \stackrel{\sim}{\longrightarrow} \VERT(T)
  \sqcup \Edge(T) \sqcup \Leg(T).
\end{equation}

\begin{proof}[Proof of Theorem \ref{legendre}]
  By definition, the sum
  \begin{multline*}
    \sum_{(T,v)\in\Gamma^{\VERT{}}_{0,n}(d)}
    \bigl[\CM_\square(T)/\Aut(T,v) \bigr] +
    \sum_{(T,e)\in\Gamma^{\Edge{}}_{0,n}(d)} \bigl[
    \CM_\square(T)/\Aut(T,e) \bigr] \\
    + \sum_{(T,i)\in\Gamma^{\Leg{}}_{0,n}(d)} \bigl[
    \CM_\square(T)/\Aut(T,i) \bigr]
  \end{multline*}
  equals
  \begin{equation}\label{sss}
    \sum_{T\in\Gamma_{0,n}(d)} \bigl[ \bigl( \bigl( \VERT(T) \sqcup
    \Edge(T) \sqcup \Leg(T) \bigr) \times \CM_\square(T)\bigr)/\Aut(T) 
    \bigr].
  \end{equation}
  By the canonical identification \eqref{zzz}, the sum \eqref{sss} may
  be rewritten as
  \begin{equation*}
    \sum_{T\in\Gamma_{0,n}(d)} \bigl[ \bigl( \bigl( \ast \sqcup
    \Flag(T) \bigr) \times \CM_\square(T)\bigr)/\Aut(T) \bigr]
  \end{equation*}
  which equals
  \begin{equation*}
    [\Mbar_{0,n}(r,d)] + \sum_{(T,f)\in\Gamma^{\Flag{}}_{0,n}(d)}
    \bigl[\CM_\square(T)/\Aut(T,f) \bigr] .
  \end{equation*}
  The theorem now follows from Lemma \ref{ccc}.
\end{proof}

\section{Composition on the Grothendieck group of mixed Hodge structures}

Let $\Cat$ be an abelian category, and let $\Gamma$ be a finite group.
The Grothendieck group $K(\Cat,\Gamma)$ is generated by elements
$[M]$, where $M=M_0\oplus M_1$ is a $\Z/2$-graded finitely generated
projective $\Gamma$-module in $\Cat$, subject to the following
relations:
\begin{enumerate}
\item $[M]=[M\oplus N]$ if $N$ is a finitely generated projective
$\Gamma$-module with $N_0\cong N_1$;
\item $[M] = [M/N]+[N]$ when $N$ is a $\Gamma$-invariant $\Z/2$-graded
submodule of $M$.
\end{enumerate}
The negative of an element $[M_0\oplus M_1]$ of $K(\Cat,\Gamma)$ is
$[M_1\oplus M_0]$. Every element of $K(\Cat,\Gamma)$ is represented by
a class $[M]$ for some $M$.

There is a generalization when $\Gamma$ is replaced by the groupoid
$\SS$: the resulting Grothendieck group decomposes as a product
\begin{equation*}
  K(\Cat,\SS) = \prod_{n=0}^\infty K(\Cat,S_n) .
\end{equation*}
If $\Cat$ has an exact symmetric tensor product $\o$, then
$K(\Cat,\SS)$ is a commutative ring, whose product is induced by the
product on $\Z/2$-graded $\SS$-modules
\begin{equation*}
  \bigl( \XX \boxtimes \YY \bigr)(n) = \bigsqcup_{m=0}^n
  \Ind_{S_m\times S_{n-m}}^{S_n} \XX(m) \o \YY(n-m) .
\end{equation*}
Here, we use the standard extension of the symmetric tensor product to
$\Z/2$-graded objects of $\Cat$, sometimes referred to as the Koszul
sign convention.

In particular, $K(\Cat,\SS)$ is an algebra over the commutative ring
\begin{equation*}
  K(\Cat) = K(\Cat,S_0) \subset K(\Cat,\SS) .
\end{equation*}
The identity element $1\in K(\Cat)$ is the class of the identity $\1$
for the tensor product on $\Cat$.

The algebra $K(\Cat,\SS)$ is filtered. The subspace $F_nK(\Cat,\SS)$
is spanned by classes $[\XX]$ such that $\XX(m)=0$ for $m<n$.

Define the composition of $\SS$-modules $\XX$ and $\YY$ in $\Cat$ by
\begin{equation*}
  \bigl( \XX\circ\YY \bigr)(n) = \bigsqcup_{i=0}^\infty \XX(i)
  \o_{S_i} \YY^{\boxtimes i}(n) .
\end{equation*}
This operation is associative:
\begin{equation*}
  (\XX\circ\YY)\circ\ZZ \cong \XX\circ(\YY\circ\ZZ) ,
\end{equation*}
and descends to an associative operation $\circ$ on $K(\Cat,\SS)\times
F_1K(\Cat_k,\SS)$ such that
\begin{equation*}
  [\XX]\circ[\YY] = [\XX\circ\YY] .
\end{equation*}

There is a functor $\XX\mapsto\delta\XX$ from $\SS$-modules to
$\SS$-modules, defined by
\begin{equation*}
  \delta\XX(n) = \Res^{S_{n+1}}_{S_n} \XX(n+1) 
\end{equation*}
satisfying the formula
\begin{equation*}
  \delta(\XX\boxtimes\YY) \cong \delta\XX\boxtimes\YY \sqcup
  \XX\boxtimes\delta\YY .
\end{equation*}
The linear operation $D$ on $K(\Cat,\SS)$, defined by
$D[\XX]=[\delta\XX]$, induces a derivation of the algebra
$K(\Cat,\Gamma)$ over $K(\Cat)$ such that
\begin{equation*}
  D : F_nK(\Cat,\SS) \to F_{n-1}K(\Cat,\SS) .  
\end{equation*}
By the same proof as for $K(\Var,\SS)$, we see that $D$ satisfies the
equation
\begin{equation*}
  D(a\circ b) = (Da\circ b) \, Db .
\end{equation*}

Let $s_n$ be the class associated to the $\SS$-space $\1$ with trivial
action of $S_n$. We have $s_0=1$, $Ds_n=s_{n-1}$. In particular,
$Ds_1=1$, and
\begin{equation*}
  s_1\circ[\XX] = [\XX]\circ s_1 = [\XX] .
\end{equation*}

Assembling the above constructions on $K(\Cat,\SS)$, we obtain the
following theorem.
\begin{theorem} \label{composition:abelian}
  If $\Cat$ is an abelian category with exact symmetric tensor
  product, the algebra $K(\Cat,\SS)$ is a complete algebra with
  composition operation over $K(\Cat)$.
\end{theorem}

Now suppose that $\Cat$ is defined over $\Q$. The representation
theory of the symmetric group in the category of vector spaces over
$\Q$ can be transferred \emph{mutatis mutandi} to $\Cat$, giving rise
in particular to a natural identification
\begin{equation*}
  K(\Cat,\SS) \cong K(\Cat)\[s_1,s_2,\dots\] .
\end{equation*}
This yields explicit formulas for $s_n\circ (ab)$ and for $s_k\circ
s_n$, which we now recall.

\begin{definition}
  A partition $\mu\vdash n$ is a decreasing sequence
  $\mu_1\ge\dots\ge\mu_\ell$ of positive integers such that
  \begin{equation*}
    n=\mu_1+\dots+\mu_\ell .
  \end{equation*}
  Denote the length $\ell$ of $\mu$ by $\ell(\mu)$.
\end{definition}

Let $s_\mu$ be the Schur polynomial associated to the partition $\mu$,
defined by the Jacobi-Trudy formula
\begin{equation*}
  s_\mu = \det\bigl(s_{\mu_i-i+j}\bigr)_{1\le i,j\le\ell(\mu)} .
\end{equation*}
The following formula for the composition operation in $K(\Cat,\SS)$
follows from the representation theory of the symmetric group (Section
I.8, \cite{Macdonald}):
\begin{equation} \label{product}
  s_n\circ(ab) = \sum_{\mu\vdash n} (s_\mu\circ a)  (s_\mu\circ b) .
\end{equation}

The Newton polynomials $p_n$ are the elements of $F_n K(\Cat,\SS)$
defined by the recursion
\begin{equation} \label{Newton}
  ns_n = p_n + s_1p_{n-1} + \dots + s_{n-1}p_1 .
\end{equation}
For example, $p_1=s_1$ and $p_2=2s_2-s_1^2$. After tensoring with
$\Q$, this formula may be inverted to give a formula for $s_n$ as a
polynomial in $\{p_1,\dots,p_n\}$. The following result is proven in
Section I.8 of \cite{Macdonald}.

\begin{lemma}
  For each $n\geq 1$, the function $a\mapsto p_n\circ a$ is a
  homomorphism of $K(\Cat,\SS)$.
\end{lemma}

The composition $s_k\circ s_n$ in $K(\Cat,\SS)$ may be calculated
using the formula
\begin{equation*}
  p_k\circ p_n = p_{kn}
\end{equation*}
for composition of Newton polynomials (Section I.8, \cite{Macdonald}).
The derivation takes the form $Dp_n=\delta_{n,1}$ on the Newton
polynomials.

In $K(\Cat,\SS)$, the formula for $\Log$ takes the simplified form
\begin{equation} \label{Log}
  \sum_{n=1}^\infty \ell_n = \sum_{k=1}^\infty \frac{\mu(k)}{k}
  \log(1+p_k),
\end{equation}
where $\mu$ denotes the M\"obius function. To see this, observe that
\begin{align*}
  \Exp \circ \Bigl( \sum_{j=1}^\infty \frac{\mu(j)}{j} \log(1+p_j)
  \Bigr) &= \exp \Bigl( \sum_{i=1}^\infty \frac{p_i}{i} \circ
  \sum_{j=1}^\infty \frac{\mu(j)}{j} \log(1+p_j) \Bigr) - 1 \\ &= \exp
  \Bigl( \sum_{i,j=1}^\infty \frac{\mu(j)}{ij} \log(1+p_{ij}) \Bigr) -
  1 = s_1 .
\end{align*}

Take $\Cat$ to be the abelian category of mixed Hodge structures
$\MH$.  The Serre characteristic
\begin{equation*}
  \Serre : K(\Var,\SS) \to K(\MH,\SS) ,
\end{equation*}
is the morphism of complete filtered algebras with composition
operation induced by taking cohomology with compact supports
\begin{equation*}
  [\XX] \mapsto \sum_i (-1)^i [H^i_c(\XX,\Q)] .
\end{equation*}
The formula for the class of the configuration space $[\FF(X,n)]$
simplifies on applying the Serre characteristic, by \eqref{Log}:
\begin{equation*}
 \Serre(\FF(X)) = \prod_{n=1}^\infty
 (1+p_n)^{\frac{1}{n}\sum_{k|n}\mu(n/k)p_k\circ \Serre(X)} .
\end{equation*}
In the special case $X=\P$, this becomes
\begin{equation*}
  \Serre(\FF(\P)) = (1+p_1)\prod_{n=1}^\infty (1+p_n)^{\frac{1}{n}
  \sum_{k|n} \mu(n/k) \L^k} ,
\end{equation*}
where $\L=\Serre(\A)$, since $p_k\circ[\P]=\L^k$ for $k\ge2$. In
particular,
\begin{equation*}
  \Serre(\FF(\P,n)) =
  \begin{cases}
    (\L+1)s_1 , & n=1 , \\
    (\L^2-\L)s_2 + \L s_1^2 , & n=2 , \\
    (\L^3-\L)s_3 , & n=3 .
  \end{cases}
\end{equation*}
We also see that
\begin{equation} \label{DFF}
  D\Serre(\FF(\P)) = \Serre(\P) \Serre(\FF(\A)) .  
\end{equation}

If $X$ is a smooth projective variety (or if $X$ has finite quotient
singularities), the Hodge numbers of $X$ are determined by the Serre
characteristic $\Serre(X)$, since the cohomology of degree $i$ is pure
of weight $i$. Thus, to calculate the Betti numbers of
$\Mbar_{0,n}(r,d)$, it suffices to calculate the Serre characteristic
$\Serre([\Mbar(r)])$, which we abbreviate to $\Serre(\Mbar(r))$.  By
Theorem \ref{legendre}, we have
\begin{equation*}
  \Serre(\Mbar(r)) = \Serre(\CM(r)) \circ ( s_1 + \Serre(\Mbar^*(r))
  \bigr) + \Serre(\P^r) \bigl( s_2 \circ\Serre(\Mbar^*(r)) -
  \Serre(\Mbar^*(r))^2 \bigr) .
\end{equation*}
It only remains to calculate $\Serre(\CM(r))$.

Since the fibration from $\Map_d(\P,\P^r)\times\FF(\P,n)$ to
$\CM_{0,n}(r,d)$ is not locally trivial, it is not possible to obtain
a closed formula for the class $[\CM_{0,n}(r,d)]$ from the formulas we
obtained above for $[\Map_d(\P,\P^r)]$ and $[\FF(\P,n)]$. On taking
the Serre characteristic, we are able to circumvent this obstruction
by use of the following theorem, whose proof uses an argument
suggested to us by B. Totaro.
\begin{theorem} \label{almost-free}
  Let $X$ be a quasi-projective variety carrying commuting actions of
  a discrete groupoid $\Gamma$ and a connected algebraic group $G$. If
  the action of $G$ on $X$ is almost-free, that is, the isotropy
  groups of the action are finite, then $\Serre(X,\Gamma) = \Serre(G)
  \, \Serre(X/G,\Gamma)$.
\end{theorem}

\begin{lemma} \label{canon}
  If $\Orbit$ is an almost-free homogeneous space for the connected
  group $G$, then for each $k\ge0$, there is a canonical isomorphism
  of mixed Hodge structures $H^k_c(\Orbit,\Q)\cong H^k_c(G,\Q)$.
\end{lemma}
\begin{proof}
  The choice of a point $x\in\Orbit$ determines an algebraic map
  $\mu_x:G\to\Orbit$, and hence a map on cohomology
  $\mu_x^*:H^k_c(\Orbit,\Q) \to H^k_c(G,\Q)$. Since $G$ is connected,
  the morphisms $\mu_x$ and $\mu_y$ are homotopic for $x,y\in\Orbit$,
  so $\mu_x^*$ is independent of $x$.

  The map $\mu^*_x$ is an isomorphism onto the subspace of
  $H^k_c(G,\Q)$ fixed by the action of the stabilizer $G_x$ of $x$.
  But $G$ is connected, so the action of $G$, and \emph{a fortiori} of
  $G_x$, on $H^k_c(G,\Q)$ is trivial.  Hence, $\mu^*_x$ is an
  isomorphism between the vector spaces $H^k_c(\Orbit,\Q)$ and
  $H^k_c(G,\Q)$. Since $\mu_x$ is an algebraic map, $\mu^*_x$ is also
  an isomorphism of mixed Hodge structures.
\end{proof}

\begin{proof}[Proof of Theorem \ref{almost-free}]
  Denote by $\pi$ the projection from $X$ to $B=X/G$. There is a
  $\Gamma$-invariant Zariski open subset $B_0\subset B$ on which
  $\pi:X_0=\pi^{-1}(X)\to B_0$ is a locally trivial fibration of
  analytical spaces. By Noetherian induction, it suffices to prove the
  theorem when $\pi:X\to B$ is a locally trivial fibration of analytic
  spaces.
  
  By M. Saito's theory of mixed Hodge modules, the Leray spectral
  sequence with compact supports
  \begin{equation*}
    E_2^{ij} = H^i_c(B,R^j\pi_!\Q_X) \Rightarrow H^{i+j}_c(X,\Q)
  \end{equation*}
  is a spectral sequence of $\Gamma$-equivariant mixed Hodge
  structures. Lemma \ref{canon} identifies the $E_2$-term with the
  tensor product of mixed Hodge structures
  \begin{equation*}
    E_2^{ij} = H^i_c(B,\C) \o H^j_c(G,\C) .
  \end{equation*}
  It follows that
  \begin{align*}
    \Serre(X,\Gamma) &= \sum_k (-1)^k [H^k_c(X,\Q)] \\ &= \sum_{i,j}
    (-1)^{i+j} [E_2^{ij}] = \Serre(B,\Gamma) \, \Serre(G) ,
  \end{align*}
  proving the theorem.
\end{proof}

We are now able to calculate the Serre characteristic of
$\CM_{0,n}(r,d)$.
\begin{theorem}
  \begin{multline*}
    \frac{\Serre(\CM(r))}{\Serre(\P^r)} =
    \frac{\Serre(\P^{r-1})}{\Serre(\P)}
    \, \frac{q\,\Serre(\FF(\P))}{1-q\L^{r+1}} \\
    + \frac{\Serre(\FF(\P)) - (1+ (\L+1)s_1 +
      (\L^2-\L) s_2 + \L s_1^2)}{\L^3-\L} .
  \end{multline*}
\end{theorem}
\begin{proof}
  The moduli space $\CM_{0,n}(r,d)$ equals the quotient
  \begin{equation*}
    \Map_d(\P,\P^r) \times_{\Aut(\P)} \FF(\P,n) ,
  \end{equation*}
  where $\Aut(\P)=\PGL(2)$ acts on both $\Map_d(\P,\P^r)$ and
  $\FF(\P,n)$ via the automorphisms of $\P$. This representation is
  compatible with the action of $S_n$.
  
  If the stability condition is satisfied, the action of $\Aut(\P)$ on
  $\Map_d(\P,\P^r)\times\FF(\P,n)$ is almost-free. Applying Theorem
  \ref{almost-free}, we see that
  \begin{equation*}
    \Serre(\CM(r,d)) =
    \frac{\Serre(\Map_d(\P,\P^r))\Serre(\FF(\P))} {\Serre(\Aut(\P))} .
  \end{equation*}
  The formula for $\Serre(\CM(r))$ is similar, except that the terms
  of degree 0 with $n=0$, $1$ and $2$ in $\Serre(\FF(\P))$ must be
  removed since they do not correspond to stable curves.  Since
  $\Serre(\Aut(\P))=\L^3-\L$, the theorem follows.
\end{proof}

%The following corollary follows by \eqref{DFF}.
%\begin{corollary}
%  \begin{equation*}
%    \Serre(\CM^*(r)) = \Serre(\P^{r-1}) \, \frac{q\,
%      \Serre(\FF(\A))}{1-q\L^{r+1}}
%  \end{equation*}
%\end{corollary}

Taking the coefficient corresponding to $n=0$, we see that
\begin{equation*}
  \Serre(\CM_{0,0}(r,d)) = \frac{\Serre(\Map_d(\P(V),\P^r))}{\L^3-\L}
  = \L^{(d-1)(r+1)} \, \frac{\Serre(\P^r)\Serre(\P^{r-1})}{\Serre(\P)}
  .
\end{equation*}
This is consistent with the isomorphism between $\CM_{0,0}(r,1)$ and
the Grassmannian $G(2,r+1)$. In fact, for all $d>0$, the spaces
$\CM_{0,0}(r,d)$ and $G(2,r+1)$ have isomorphic Chow rings
\cite{rahul:chow}.

It is interesting to note that the Serre characteristic
$\Serre(\CM(r,d))$ stabilizes as $r\to\infty$:
\begin{multline*}
  \lim_{r\to\infty} \Serre(\CM(r)) = \frac{q\,
    \Serre(\FF(\P))}{(1-\L)(1-\L^2)} \\
  - \frac{\Serre(\FF(\P)) - ( 1 + (\L+1)s_1 + (\L^2-\L) s_2 + \L
    s_1^2)}{\L(1-\L)(1-\L^2)} .
\end{multline*}
In the limit $r\to\infty$, the Serre characteristic $\Serre(\CM(r,d))$
vanishes if $d>1$.

\appendix

\section{Calculation of $\Serre(\Mbar_{0,0}(r,3))$}

We present here an explicit calculation of the Serre
characteristic of $\Mbar_{0,0}(r,3)$. There are four strata,
corresponding to the four trees in $\Gamma_{0,0}(3)$.

\medskip

\noindent \setlength{\unitlength}{0.0075in}   \begin{picture}(20,35)(0,0)
    \thicklines
    \put(10,10){\circle*{10}}
    \put(7,20){3}
  \end{picture}

The open stratum is $\CM_{0,0}(r,3)$, and has Serre characteristic
\begin{equation*}
  \frac{\L^{2(r+1)} \Serre(\P^r)\Serre(\P^{r-1})}{\Serre(\P)} .
\end{equation*}

\bigskip

\noindent \setlength{\unitlength}{0.0075in}   \begin{picture}(100,35)(0,0)
    \thicklines
    \put(10,10){\circle*{10}}
    \put(90,10){\circle*{10}}
    \put(10,10){\line(1,0){80}}
    \put(7,20){2}
    \put(87,20){1}
  \end{picture}

This stratum is a fibre bundle over $\CM_{0,1}(r,2)$ with fibre
$\CM^*_{0,0}(r,1)$, and has Serre characteristic
\begin{equation*}
 \Serre(\CM_{0,1}(r,2)) \circ \Serre(\CM^*_{0,0}(r,1)) = \frac{(\L+1)
 \L^{r+1} \Serre(\P^r)\Serre(\P^{r-1})^2}{\Serre(\P)} .
\end{equation*}

\bigskip

\noindent \setlength{\unitlength}{0.0075in}   \begin{picture}(180,35)(0,0)
    \thicklines
    \put(10,10){\circle*{10}}
    \put(90,10){\circle*{10}}
    \put(170,10){\circle*{10}}
    \put(10,10){\line(1,0){80}}
    \put(90,10){\line(1,0){80}}
    \put(7,20){1}
    \put(87,20){1}
    \put(167,20){1}
  \end{picture}

This stratum is the quotient of a fibre bundle over $\CM_{0,2}(r,1)$
with fibre $\CM^*_{0,0}(r,1)^2$ by the automorphism group $\SS_2$ of the
tree, and has Serre characteristic
\begin{align*}
  \Serre(\CM_{0,2}(r,1)) \circ \Serre(\CM^*_{0,0}(r,1)) &= ( (\L^2-\L)
  s_2\circ[r] + \L\,s_1^2\circ[r] ) \,
  \frac{\Serre(\P^r)\Serre(\P^{r-1})}{\Serre(\P)} \\
  &= \frac{\L^2 \bigl( \Serre(\P^r) + \Serre(\P^{r-2}) \bigr)
    \Serre(\P^r)\Serre(\P^{r-1})^2}{\Serre(\P)^2} .
\end{align*}

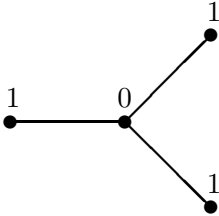
\setlength{\unitlength}{0.0075in}   \begin{picture}(155,150)(-5,0)
    \thicklines
    \put(85,65){\circle*{10}}
    \put(5,65){\circle*{10}}
    \put(145,125){\circle*{10}}
    \put(145,5){\circle*{10}}
    \put(5,65){\line(1,0){ 80}}
    \put(85,65){\line(1,1){ 60}}
    \put(85,65){\line(1,-1){ 60}}
    \put(2,75){1}
    \put(80,75){0}
    \put(142,135){1}
    \put(142,15){1}
  \end{picture}
 The last stratum is
the quotient of a fibre bundle over $\CM_{0,3}(r,0)$ with fibre
$\CM^*_{0,0}(r,1)^3$ by the automorphism group $\SS_3$, and has Serre
characteristic
\begin{align*}
  \Serre(\CM_{0,3}(r,0)) \circ \Serre(\CM^*_{0,0}(r,1)) &=
  \frac{\Serre(\P^r)\Serre(\P^{r-1})s_3\circ\Serre(\P^{r-1})}{\Serre(\P)} \\
  &=
  \frac{\Serre(\P^{r+1})\Serre(\P^r)^2\Serre(\P^{r-1})}{\Serre(\P^2)
    \Serre(\P)}
\end{align*}

When $r\to\infty$, the only two strata which contribute to
$\Mbar_{0,0}(r,3)$ are the last two:
\begin{equation*}
  \Serre(\Mbar_{0,0}(\infty,3)) =
  \frac{1+\L^2-2\L^5}{(1-\L)^2(1-\L^2)^2(1-\L^3)} .
\end{equation*}
This gives an upper bound for the Betti numbers of $\Mbar_{0,0}(r,3)$
for finite $r$.

Here is a table of Serre polynomials for small values of $r$ and for
$r=\infty$; we give the vector of Betti numbers
$(b_0,b_2,b_4,\dots)$. Note that Poincar\'e duality holds for finite
$r$, as it must.

\begin{equation*}
\begin{tabular}{|c|l|} \hline
  $r$ & $\Serre(\Mbar_{0,0}(r,3))$ \\ \hline $1$ & $(1,1,2,1,1)$ \\
  $2$ & $(1,2,5,7,9,7,5,2,1)$ \\ $3$ &
  $(1,2,6,10,17,20,24,20,17,10,6,2,1)$ \\ $\infty$ &
  $(1,2,6,11,21,32,51,71,101,133,177,223,284,\dots)$ \\ \hline
\end{tabular}\end{equation*}

\end{document}